\input amstex
 
\documentstyle{amsppt}\magnification=1200
\NoBlackBoxes
 
\magnification\magstep1
\pagewidth{32.622pc}
\pageheight{44.8pc}
\define\oz{\overline{\omega}}
\define\Oz{\overline{\Omega}}
\define\om{\overline{\omega\omega}}
\define\Om{\overline{\Omega\Omega}}
\define\ind{\operatorname{ind}}
\define\Ind{\operatorname{Ind}}

\define\ca{\Cal C_{c}^{\infty}}
 
\topmatter
\title Exceptional $\Theta$-correspondences I 
\endtitle
\author Kay Magaard and Gordan Savin \endauthor
\address  Department of Mathematics, Wayne State University,  
Detroit, MI 48202 \endaddress
\email kaym\char'100math.wayne.edu \endemail
\address  Department of Mathematics, University of Utah, Salt Lake City, UT
84112\endaddress
\email savin\char'100math.utah.edu\endemail
\thanks Research at MSRI is
   supported in part by NSF grant DMS-9022140\endthanks
\abstract Let $G$ be a split simply laced group defined over 
a $p$-adic field $F$. In this paper we study the restriction of the 
minimal representation
of $G$ to various dual pairs in $G$. For example, the restriction 
of the minimal representation of $E_7$ to the dual pair 
$G_2 \times Sp(6)$ gives the non-endoscopic Langlands lift of 
irreducible representations of $G_2$ to $Sp(6)$. 
\endabstract 
\endtopmatter
 
\document

\head Introduction \endhead 

Let $F$ be a $p$-adic field  
and  $G$ be a split simple group of type $D_n$ 
or $E_n$, ($n=6,7,8$) over $F$. If the type is $D_n$ then  
$G$ is  $SO(2n)$. In \cite{KS} the minimal representation $V$ 
of $G$ was introduced and in \cite{S} some further properties of 
it were studied. It is an analogue of the Weil representation of 
$Sp(2n)$ and in some ways it is better behaved then the Weil 
representation. Therefore it makes sense to ask if one can 
use it to obtain dual pair correspondences.

We have two types of dual pairs in $G$. In $SO(2n+2)$ 
we have a dual pair 
$$SO(2n-1)\times SO(3).$$
If $G$ is exceptional then we have a dual pair 
$$
G_2(F)\times H_D 
$$
with $H_D$ of type $A_2$, $C_3$ or $F_4$, respectively. 

Let $A\times B$ be a dual pair in $G$ and $X$ an irreducible 
representation of $A$. We say that an irreducible representation
$Y$ of $B$ is a $\Theta$-lift of $X$ if $X\otimes Y$ is a quotient
of $V$ (see \cite{H1}). Let $\Theta(X)$ be the set of all such $Y$. 

In this paper we compute Jacquet functors $V_U$ of $V$ where $U$ 
is a unipotent radical of a maximal parabolic subgroup of
$SO(3)$ or $H_D$, with $U$ commutative, or  
a unipotent radical of a maximal parabolic subgroup of $G_2$ 
such that $U$ is a Heisenberg group. 
The results are reminiscent of those of Kudla 
\cite{Ku} and Rallis \cite{Ra} in the classical case. 
We then use the knowledge of Jacquet functors to 
obtain dual pair correspondences. For example, we compute  
$\Theta$-lifts of tempered spherical representations of $SO(3)$,
$A_2$, $C_3$ and $G_2$. The lifts from $A_2$ to $G_2$ and 
from $G_2$ to $C_3$ obtained by restricting the 
minimal representation of $E_6$ and $E_7$ respectively, are Langlands 
correspondences. 

In a sequel to this paper we plan to compute Jacquet functors  
$V_U$ for unipotent radicals of remaining maximal parabolic 
subgroups. Once this is done, it will be possible to compute 
lifts of most Langlands quotients and check the conjectural 
answer for the exceptional correspondences as given by B. Gross
\cite{G2}.
Local computations are, in a way, a preparation for global 
correspondences (i.e. correspondences of automorphic forms). So 
it is worth mentioning that 
in a forthcoming work, D. Ginzburg, S. Rallis and D. Soudry 
are studying a global variant of exceptional correspondences. 
Another possible global application has recently been initiated by 
B. Gross in connection to a realization of a $G_2$-motive 
\cite{G1}. 

The paper is organized as follows. In Section 1 we recall the definition and 
few properties of the minimal representation from \cite{S}. In Section 2 
we study the dual pair $SO(3)\times SO(2n-1)$. It is the simpliest case 
and as such it is a good introduction to exceptional 
dual pairs which form a more interesting part of this work. In 
Sections 3 and 4 we study the dual pair $G_2\times A_2$.  We compute 
$\Theta$-lifts of spherical tempered representations of $A_2$. 
In Section 5 we study the dual pair $G_2\times C_3$.  We compute 
$\Theta$-lifts of spherical tempered representations of $C_3$.
In Section 6 we recall some facts about simple Jordan  algebras 
or rank $3$. Exceptional groups can be described nicely in terms of
Jordan algebras. We use this description in Section 7 to compute $V_U$ 
where $U$ is the unipotent radical of the "Heisenberg" maximal parabolic
subgroup of $G_2$ in all three cases. We finish the paper by 
computing $\Theta$-lifts of spherical tempered representations of 
$G_2$ in all three cases. 

\head 1. Minimal representation \endhead 

Let $G$ be a simple, split, group of type $A_{2n-1}$, $D_n$ or $E_n$.  
Let $^LG$ be the complex $L$-group of $G$ (see \cite{B}). It is well known that
conjugacy classes of unipotent elements in $^LG$ correspond to 
conjugacy classes of homomorphisms
$$
\varphi : SL_2(\bold C) \rightarrow\, ^LG.
$$
Assume now that $\varphi$ corresponds to the 
subregular unipotent orbit. Let $q$ be the order of the residual 
field of $F$. Let  
$$
s = \varphi 
\left( \smallmatrix 
q^{\frac12} & 0 \\
0 & q^{-\frac12}
\endmatrix \right).
$$
Then $V$ is the spherical representation of $G$ with the 
Satake parameter $s$ (see \cite{Ca}). 

Let $\frak g$ be the Lie algebra of $G$. Let 
$<\cdot , \cdot >$ be the Killing form on $\frak g$.  
Fix $\psi : F \to \Bbb C^{\times}$, a non-trivial character. 
Let $\ca(\frak g)$ denote the space of locally constant, 
compactly supported functions on $\frak g$. Define the Fourier 
transform on $\ca(\frak g)$ by 
$$
\hat f(y)=\int_{\frak g} f(x) \psi(<x,y>) dx .
$$
Let $\Cal O_{\min}$ be the unique minimal nilpotent $G$-orbit in 
$\frak g$ and $\mu_{\Cal O_{\min}}$ a $G$-invariant measure on 
$\Cal O_{\min}$ normalized as in \cite{MW}. 
Let $\frak g_{\Bbb Z}$ be the $\Bbb Z$-span of a Chevaley basis 
of $\frak g$. Let $\Cal R$ be the ring of integers of $F$ and 
$\varpi$ the uniformizing element. Let 
$$
\frak g _n = \frak  g_{\Bbb Z}\otimes_{\Bbb Z} \varpi^n \Cal R.
$$  
In \cite{S} we have shown that there exists $m>0$ such that 
$$
tr\int_{\frak g}f(x)\pi(\exp x) dx = \int \hat f \mu_{\Cal O_{\min}}
$$
for any $f \in \ca(\frak g)$ supported on $\frak g_m$ and such that 
$\hat f(0)=0$. 

This implies "smallness" of
the $N$-spectrum of $V$, where $N$ is a unipotent radical of
a maximal parabolic subgroup of $G$. 
More precisely, let  
$\Delta$ be the Dynkin diagram of $G$. Mark  
the diagram $\Delta$ as follows. Attach $0$ to the unique 
branch vertex (or to the middle vertex of $\Delta$ if  
the type of $G$ is $A_{2n-1}$) and $2$ to all other. This marking corresponds 
to the subregular nilpotent orbit \cite{D}.  Let $P=MN$ be 
a maximal parabolic subgroup of $G$. Let $\Delta_M$ be the 
Dynkin diagram of $M$. Assume that we are in the 
following favorable situation: 
\roster 
\item $N$ is a commutative group.  
\item The marking of $\Delta$ corresponding to the subregular
nilpotent orbit of $G$ restricts to the marking of $\Delta_M$ 
corresponding to the subregular nilpotent orbit of $M$. 
\endroster 
The possible cases are given by the following table:

$$
\matrix 
G & M & N \\
D_{n+1} & D_{n} & F^{2n} \\
E_6 & D_5 & F^{16} \\
E_7 & E_6 & F^{27}
\endmatrix
$$
Here $F^{2n}$ is the standard representation of $SO(2n)$, 
$F^{16}$ is a spin-representation of $Spin_{10}$ and $F^{27}$
is isomorphic to the exceptional Jordan algebra. We say that a 
point in $N$ is singular if it is a highest weight vector for 
a Borel subgroup of $M$. Let $\omega$ and 
$\overline{\omega}$ be the sets of singular vectors in $N$ and 
$\overline{N}$. Note that $\omega$ is the smallest non-trivial 
$M$-orbit in $N$. If $G=SO(2n)$ then $\omega$ is the null-cone 
of the invariant quadratic form (with $0$ excluded). 

\proclaim{Theorem 1.1} The minimal representation $V$ 
of $G$ 
has a $P$-invariant filtration 
$$
0 \to {\Cal C}_{c}^{\infty}(\oz) \to V 
\to V_N \to 0. 
$$ 
Here $\Cal C_{c}^{\infty}(\oz)$ denotes the space of locally
constant, compactly supported functions on $\oz$ and $V_N$ is 
the space of $N$-coinvariants of $V$ (Jacquet functor).
\roster 
\item Let $f\in \Cal C_{c}^{\infty}(\oz)$. The action of $P$ is given by  
$$\
\pi(n)f(x) = \psi(<x,n>)f(x), \qquad n\in N  
$$
and
$$
\pi(m)f(x) = |\det(m)|^{\frac{s}{d}}f(m^{-1}x), \qquad m\in M. 
$$
\item 
$$
V_N \cong V_M \otimes |\det|^{\frac{t}{d}} + 
 |\det|^{\frac{s}{d}}
$$ 
where $V_M$ is the minimal representation of $M$ (center acting 
trivially). 
\endroster 
In the above formulas 
$<\cdot , \cdot >$ is the $F$-valued pairing between 
$N$ and $\overline N$ induced by the Killing form on $\frak g$, 
$\det$ is determinant of 
the representation of $M$ on $N$, $d$ is the dimension of $N$ and  
the values of $s$ and $t$ are given in the following table.  
$$
\matrix 
G & s & t \\
D_{n+1} & n-1 & 1 \\ 
E_6 & 4 & 2 \\ 
E_7 & 6 & 3 
\endmatrix 
$$
\endproclaim 

\demo{Proof} 
This is just Theorem 6.5 in \cite{S} if $G$ is $E_7$. 
The other two cases also satisfy conditions of 
Proposition 4.1 in \cite{S}. Hence
the proof carries over with
no changes. The proof given in \cite{S}, however, is 
valid only when the residual characterisitic is odd. 
It remains to discuss the case when the residual characteristic is $2$. 

Let $x\in \overline N$, and define $\psi_x(n)=\psi(<x,n>)$. Let 
$V_{N,\psi_x}$ be the quotient of $V$ by the space spanned by the 
elements $\{\pi(n)v-\psi_x(n)v | n\in N, v \in V\}$. The key point in 
the proof of Theorem 6.5 in \cite{S} is to show that 
$$
V_{N,\psi_x}=0 
$$
for $x\neq 0$ and not in $\overline \omega$, i.e. the $N$-spectrum of $V$
is concentrated on the closure of the smallest $M$-orbit in $\overline N$. 
This follows from the 
work of Moeglin and Waldspurger \cite{MW},  
but only if the residual characteristic of $F$ is odd. 

Let $k$ be a number field and $\Bbb A$ its ring of adeles. 
Ginzburg, Rallis and Soudry \cite{GR} have constructed a square 
integrable automorphic form $\Pi=\otimes_{v}\Pi_{v}$ 
on $G_{\Bbb A}$ isomorphic to $V$ at every finite 
place. But arguing exactly as Howe (\cite{H2} 
Lemma 2.4), one shows that if the $N$-spectrum is concentrated on 
the closure of $\overline \omega$  at one place $v$, it has to 
be concentrated on the closure of $\overline \omega$  at all 
places. This completes the proof of the theorem.  

\enddemo

\head 2. Dual pair $SO(2n-1) \times SO(3)$ \endhead 

Let $G=SO(2n+2)$. 
 Let $e,h,f$ be the standard basis 
for $sl(2)$, the Lie algebra of $SO(3)=PGL_2(F)$. Let $P=MN$ be the maximal parabolic
of $G$ as in Section 1. Let $\frak g$ be the Lie algebra of $G$ and let 
$$
\frak g = \overline{\frak n} + \frak m + \frak  n
$$ 
be the corresponding decomposition of $\frak g$. The embedding $SO(3) 
\subseteq SO(2n+2)$ corresponds to $sl(2)\subseteq \frak  g$  such that 
\roster 
\item 
$\overline{\frak n}=\{ x\in \frak g~ |~ [h,x]=-2x\}. $ 
\item 
$\frak m=\{ x\in \frak g~ |~ [h,x]=0\}. $ 
\item 
${\frak n}=\{ x\in \frak g~ |~ [h,x]=2x\}. $ 
\endroster 
Note that $e\in \frak n$, $f \in \overline{\frak n}$ and 
$SO(2n-1)$ is the centralizer of $e$ in $M$.

Let $Q=LU =P\cap SO(3)$ be a Borel subgroup of $SO(3)$. 

\proclaim{Proposition 2.1} 
Let $\pmatrix  a & \\ & b \endpmatrix \in L \subset PGL_2(F)=SO(3)$. 
\roster 
\item 
$$
\pi(\pmatrix  a & \\ & b \endpmatrix)f(x) = |\frac{a}{b}|^{n-1} f(\frac{a}{b} x), \qquad f\in 
\ca (\oz). 
$$
\item 
The eigenvalues of $\pi(\pmatrix  a & \\ & b \endpmatrix)$ on 
$V_N$ are $|\frac{a}{b}|$ and 
$|\frac{a}{b}|^{n-1}$. 
\endroster 
\endproclaim 

\demo{Proof} This is just Theorem 1.1. 
\enddemo 

Let $\chi$ be a multiplicative character of $F$. Let $\rho\chi$ 
denote the character of $L$ defined by 
$$
\rho\chi(\pmatrix  a & \\ & b \endpmatrix) 
= |\frac{a}{b}|^{\frac{1}{2}}\chi(\frac{a}{b}). 
$$
Let $\tau_{\chi}= \Ind_Q^{SO(3)}\rho\chi$. If $\chi$ is unitary than 
$\tau_{\chi}$ is an irreducible tempered spherical representation of 
$SO(3)$. 
Let $\sigma \in \Theta(\tau_{\chi})$. 
By Frobenius reciprocity
$$
Hom_{SO(2n-1)\times SO(3)}(V, \sigma\otimes \tau_{\chi})= 
Hom_{SO(2n-1) \times L}(V_U, \sigma\otimes \rho\chi).
$$
Therefore, $\sigma\in \Theta(\tau_{\chi})$ is a quotient of $V$ if and 
only if $\sigma\otimes \rho\chi$
is a quotient of $V_U$. Since
$$
0 \to {\Cal C}_{c}^{\infty}(\oz)_U \to V_U  
\to V_N \to 0 
$$ 
we need to understand $\ca(\oz)_U$. 
Let $NN$ and $\overline{NN}$ be the 
complements of $\overline U$ and $U$ in $N$ and $\overline N$. 
Let $\overline \omega$, as before, be the 
set of singular vectors in $\overline N$. Put 
$$
\om =\overline{\omega} \cap \overline{NN}. 
$$

\proclaim{Lemma 2.2}

\endproclaim 
$$
\ca(\overline\omega)_U=\ca(\om).
$$
\demo{Proof} 
Let us recall few known facts about Jacquet functors. Let $E$ be a $U$-module. Then 
$E_U= E / E(U)$ where $E(U)$ can be defined either as the space spanned by the 
elements $\{\pi(u)v - v | u\in U, v\in E\}$ or the space of all $v$ such that 
$$
\int_{U_c}\pi(u)v du = 0
$$
for some open compact subgroup $U_c \subset U$ depending on $v$ (2.33 \cite{BZ}). 

Obviously, $\ca(\om)$ is a quotient of $\ca(\overline\omega)$ and by Theorem 1.1 (1),  $U$ acts 
trivially on $\ca(\om)$. 
Let $f\in \ca(\overline\omega)$ such that $f|_{\om}=0$. To prove the lemma, we need to 
find an open compact subgroup $U_c$ such that 
$$
\int_{U_c}\psi(<x,u>)f(x)du=0
$$
for all $x\in \overline\omega$. Let $x$ be such that $f(x)\neq 0$. Since $x$ is not 
in $\overline{NN}$, there exists an open compact neighborhood $\Cal O$ of $x$ and 
an open compact subgroup $U_{\Cal O}$ such that $\psi(<y,u>)$ is a non-trivial 
character of $U_{\Cal O}$ for any $y\in \Cal O$. Since the support of $f$ is 
compact, a finite collection of $\Cal O$ covers the support of $f$. The union of the 
corresponding $U_{\Cal O}$ is the desired $U_c$. The lemma follows.
\enddemo

We can, 
therefore, summarize the situation with the following 
proposition. 

\proclaim{Proposition 2.3} 
$V_U$ has a filtration with two quotients:
$$
\ca(\om), \, \text { and } V_N.
$$
As $SO(2n-1)\times L$-modules: 
\roster 
\item the action on $\ca(\om)$ is geometric, twisted with 
$
|\frac{a}{b}|^{n-1} 
$ for $ 
\pmatrix a & \\ & b \endpmatrix \in L$. 
\item $$V_N\cong V_M \otimes  |\frac{a}{b}| + 
1\otimes |\frac{a}{b}|^{n-1}
$$ 
where $V_M$ is the minimal representation of $M$ (center 
acting trivially). 
\endroster 
\endproclaim

Note that $\overline{NN} = F^{2n-1}$ and $\om$ is the null-cone of the $SO(2n-1)$-invariant 
quadratic form (with $0$ excluded). Let $\Cal C^{\infty}(\om)$ be the space of locally constant
functions on $\om$. 
We can define degenerate principal series 
representations $\sigma_{\chi}$ by
$$ 
\sigma_{\chi}=\{ f\in \Cal C^{\infty}(\om) \mid f(cx)=\chi(c)
|c|^{\frac{3}{2}-n}f(x)\}. 
$$
Analogously, $\sigma_{\chi}$ can be defined as a quotient of 
$\ca(\om)$ consisting of $\overline f$ such that 
$$
\overline f(cx)=\chi(c) |c|^{\frac{3}{2}-n}\overline f(x).
$$
If $\chi$ is unramified and unitary than $\sigma_{\chi}$ is an 
irreducible unitarizable spherical representation by a result of 
Tadi\'c \cite{T2}, Theorem 9.2. We are now ready to 
state and prove the main result of this section.

\proclaim{Proposition 2.4} Let $\chi$ be an unramified, unitary 
multiplicative character. Then 
$$
\Theta(\tau_{\chi})= \{\sigma_{\chi}\}.
$$ 
\endproclaim

\demo{Proof} Let $\sigma \in \Theta(\tau_{\chi})$. Then 
$\sigma \times \rho\chi$ is a quotient of $V_U$. Since 
$\rho\chi$ is different from $|\cdot|^{n-1}$ and $|\cdot|$
if $\chi$ is unitary, it follows that $\sigma\times \rho
\chi$ is quotient of $\ca(\om)$. Therefore $\sigma$ must be 
isomorphic to $\sigma_{\chi}$. To finish the theorem, 
we have to show that $\sigma_{\chi} \in \Theta(\tau_{\chi})$. By 
Frobenius reciprocity, it suffices to show that 
$\sigma_{\chi} \times \rho\chi$ is a quotient of $V_U$. We know
that it is a quotient of $\ca(\om)$.  We need the  
following lemma. 

\proclaim{Lemma 2.5} Let $H$ be a group and 
$$
0\to V_1 \to V_2 \to V_3 \to 0
$$ 
a sequence of $H$-modules. Let $W$ be a $H$-invariant quotient of 
$V_1$. Assume that there exists a $H$-invariant linear map $T$ on 
$V_2$ reducing to scalar multiplications by $\lambda$ and $\mu$ on
$W$ and $V_3$. If $\lambda \neq \mu$ then $W$ is a quotient of $V_3$.
\endproclaim 

\demo{Proof} Let $V_0$ be such that 
$$
0\to V_0 \to V_1 \to W \to 0. 
$$ 
Then 
$$
0\to W \to V_2/V_0 \to V_3 \to 0
$$ 
Since the eigenvalues of $T$ are different, the last sequence splits and 
the lemma follows. 
\enddemo 

Let $\varpi$ be the uniformizing element of $F$. Put 
$$
T =  \pi  
\pmatrix \varpi & \\ & 1 \endpmatrix. 
$$ 
If $\chi$ is unitary, $T$ satisfies the conditions of Lemma, 
hence $\sigma_{\chi}\times \rho\chi$ is a quotient of $V_U$ and the 
theorem is proved. 
\enddemo

\head 3. Group $E_6$ \endhead 

Let $J$ be the 27-dimensional exceptional Jordan algebra and 
$(x,y,z)$ be the Dickson trilinear form on $J$ (see \cite{J2}). 
Let $G_S$ be the isometry 
group of the Dickson form. It is a split, simply-connected  
groups of type $E_6$ \cite{A1}.

Let $\frak g_S$ be the Lie algebra of $G_S$. 
Then 
$$
\frak g_S = \bar \frak n + \frak m_S + \frak n
$$
where $[\frak m_S, \frak m_S]= so(10)$ and $\bar \frak n$, $\frak n$ are 
two spin-representations. 

Let $P_S=M_S N$  
be the corresponding maximal parabolic subgroup in $G_S$. 
Then $[M_S,M_S]=Spin_{10}$ and $Spin_{10}$ can be defined as 
the subgroup of $G_S$ stabilizing the $10$-dimensional 
subspace $J_{10}$ of $J$ consisting of matrices in $J$ such that
the coefficients in the third row and the third column are zero and
fixing the vector 
$$
d= \pmatrix  
0 & 0 & 0 \\
0 & 0 & 0 \\ 
0 & 0 & 1 
\endpmatrix . 
$$
In particular, $Q_{d}(x,y)=(x,y,d)$ is the invariant quadratic 
form on $J_{10}$ and $N$ (one of the spin-modules) can be identified with 
the space in $J$ consisting of matrices 
$$
\pmatrix 
0 & 0 & y \\
0 & 0 & z \\
\overline{y} & \overline{z} & 0 
\endpmatrix . 
$$
The group $G_S$ has a unique smallest, non-trivial orbit $\Omega$
on $J$.  It is the orbit of $d$ and its dimension 
is $17$. On the other hand, $M_S$ has two non-trivial orbits on 
$N$ and the smaller (call it $\omega$) can be identified with 
$\Omega \cap N$ \cite{A1}.

\head 4. Dual pair $G_2(F) \times PGL_3(F)$ \endhead 
 
Let $D$ be a nondegenerate subalgebra of the Octonions and 
let $J_D$ denote the subalgebra of $J$ obtained by restricting the 
off-diagonal entries of the matrices of $J$ to $D$. When 
$\dim(D)= 1$, i.e. $D = <1_{\Bbb O}> = F$ then $J_F$ is the algebra of 
symmetric $3\times 3$ matrices over $F$ and $C_{G_S}(J_F)= G_2(F)$, the split
group of type $G_2$.
The centralizer of $G_2(F)$ in $G_S$ is $SL_3(F)$ \cite{A1}. 

Let $G$, $G_S\subset G$, be a reductive group of type $E_6$ appearing
as a Levi factor of a maximal parabolic subgroup of a split, 
simply-connected group of type $E_7$. Then $C_G(G_2(F))=GL_3(F)$.
The action of $G_S$ on $J$ can be extended to $G$. 
Indeed, 
the inclusion $GL_3 \subset G$ induces an isomorphism 
$GL_3/[GL_3,GL_3]\cong G/[G,G]$. Since $[G,G]=G_S$, 
it suffices to extend the action to $GL_3 \subset SL_3$.  
This is done by the formula 
$$
gxg^t
$$
where $x\in J$, $g\in GL_3(F)$ and $g^t$ is the transpose of $g$. 

Let $\frak g$ be the Lie algebra of $G$. Then, as in Section 3, 
$$
\frak g =  \bar \frak n + \frak m + \frak n
$$
with $\frak m_S \subset \frak m$. Let $P=MN$ be the corresponding maximal
parabolic subgroup of $G$. Then $Q=LU=P\cap GL_3(F)$ is a maximal 
parabolic subgroup of $GL_3$ with $L=GL_2(F)\times GL_1(F)$. One 
can write 
$$
gl(3) = \bar \frak u + \frak l + \frak u 
$$
where $\frak l \subset \frak m$, $\overline{\frak u} \subset \overline{\frak n}$ and 
$\frak u \subset \frak n$. 

\proclaim{Proposition 4.1} Let $V$ be the minimal representation of 
$G$, the center acting trivially. 
Identify $\overline{N}$ with the set of pairs of 
Octonions $(y,z)$. Let $f\in \ca(\overline\omega)$. Then: 
\roster 
\item 
$$
\pi(g)f((y,z)) = f((g^{-1}y, g^{-1}z)), \qquad g\in G_2(F).  
$$
\item 
$$
\pi(l_2 \times l_1)f((y,z)) = \frac{|\det l_2|^2}{|l_1|^4}
f(l_1^{-1}(y,z)l_2), \qquad l_2\times l_1 \in GL_2(F)\times GL_1(F)  
$$
\endroster 
where $\det$ denotes the usual determinant of $2\times 2$ 
matrices. 
\endproclaim 

\demo{Proof} This is just Theorem 1.1. 
\enddemo 
Let $Y$ be an irreducible representations of $GL_3(F)$ 
and $\tilde Y$ an irreducible representation of $G_2(F)$ such that
$\tilde Y\in \Theta(Y)$. Assume that $Y$ is a submodule of 
$\Ind_{Q}^{GL_3}X$ for some irreducible representation $X$ of $L$. 
By Frobenius reciprocity
$$
Hom_{G_2\times GL_3}(V, \tilde Y\otimes \Ind X)= 
Hom_{G_2 \times L}(V_U, \tilde Y\otimes X).
$$
Therefore, $\tilde Y\otimes Y$ is a quotient of $V$ only if $\tilde Y\otimes X$
is a quotient of $V_U$. Since
$$
0 \to {\Cal C}_{c}^{\infty}(\oz)_U \to V_U  
\to V_N \to 0 
$$ 
we need to understand $\ca(\oz)_U$. Let $NN$ and $\overline{NN}$ be the 
orthogonal complements of $\overline{U}$ and $U$ in $N$ and 
$\overline{N}$. 
Since $U$   
can be identified with 
the space in $J$ consisting of matrices 
$$
\pmatrix 
0 & 0 & y \\
0 & 0 & z \\
y & z & 0 
\endpmatrix
$$ 
where $y,z \in F$,  
it follows that 
$$
NN = \{ 
\pmatrix 
0 & 0 & y \\
0 & 0 & z \\
\overline{y} & \overline{z} & 0 
\endpmatrix 
|~ \overline{y} =-y \text { and } \overline{z} =-z \}.
$$  
Let $\om =\oz \cap \overline{NN}$. As in Lemma 2.2  
$
\ca(\oz)_U = \ca(\om ), 
$
and we have to understand the structure of $G_2(F) \times L$ 
orbits on $\om$ which is, of course, the 
same as the structure of orbits on $\omega\omega=\omega \cap NN$. 

\proclaim{Proposition 4.2} 
\roster 
\item $\om = \{(y,z)\neq (0,0)| 
 \overline{y} =-y,  \overline{z} =-z \text { and } y^2=z^2=yz=0 \}$.
\item Let $AA$ and $BB$ be the subsets  
of $\om$ consisting 
all pairs $(y,z)$ such that the space $Fy+Fz$ has dimension $2$ and $1$  
respectively. 
Then $AA$ and $BB$ are $G_2(F)\times GL_2(F)$-orbits.
\endroster 
\endproclaim 

\demo{Proof} Let us prove the corresponding statement for 
$\omega\omega$. Let   
$$
n= 
\pmatrix 
0 & 0 & y \\
0 & 0 & z \\
-y & -z & 0 
\endpmatrix \in \omega\omega .
$$
Since $n$ is a traceless matrix in $J$ it will be singular iff 
$n^2=0$. But this is equivalent to  
$y^2 =z^2= yz = 0$. The first part of the proposition is proved.

We go on to observe that $G_2$ has three orbits on the set of spaces 
of traceless Octonions with the property that the Octonion
multiplication is trivial. These are characterized by their dimension;
the possible choices being $0,1,2$.
The stabilizers of the nontrivial spaces are the maximal parabolics
of $G_2$. 

Let $x$ be a traceless Octonion such that $x^2=0$. Let
$P_1$ be the maximal parabolic subgroup of $G_2(F)$ stabilizing
the line $Fx$. The Levi factor of $P_1$ is "spanned" by a long root. Consider
$$
B=\{ (ax,0)|a\in F \text { and } a\neq 0\}. 
$$
Let $QQ$ be the maximal parabolic subgroup of $GL_2(F)$ 
stabilizing $B$. 
Then $P_1 \times QQ$ acts transitively on $B$ and 
$$
BB=  (G_2(F)\times GL_2(F)) \times_{(P_1 \times QQ)} B.
$$
Let $y$ and $z$ be two traceless and linearly 
independent Octonions such that $y^2=z^2=yz=0$. Let  
$P_2$ be the maximal parabolic subgroup of $G_2(F)$ stabilizing
the space $Fy+Fz$. The Levi factor of $P_2$ is "spanned" by a short 
root. Consider
$$
A=\{ (ay+bz,cy+zd)|a,b,c,d\in F \text { and } ad-bc\neq 0\}. 
$$ 
Then $P_2 \times GL(2)$ acts transitively on $A$ and 
since $G_2(F)$ acts transitively on the set of all 
two-dimensional spaces of traceless Octonions with 
trivial multiplication,
$$
AA =  G_2(F) \times_{P_2} A.
$$
The proposition is proved.  
\enddemo 

We can now summarize the structure of $V_U$ as a $G_2(F) \times 
GL_2(F)$-module in the following theorem (compare \cite{Ku}). 

\proclaim{Theorem 4.3}
$V_U$ has a filtration with three quotients:
$$
\ca(AA),\,\, \ca(BB),\, \text { and } V_N.
$$
Moreover 
\roster 
\item 
$$
\ca(AA) = \ind_{P_2}^{G_2}(\ca(A))\otimes |\det|^2  
$$ 
\item  
$$
\ca(BB)=\ind_{P_1\times QQ}^{G_2\times GL_2}(\ca(B))\otimes |\det|^2
$$
\item 
$$
V_N\cong  V_M \otimes |\det| + 1\otimes |\det|^2 
$$
\endroster 
as $G_2(F) \times GL_2(F)$-modules. 
\endproclaim
 
We are now ready to state and prove the main theorem. 

\proclaim{Theorem 4.4} Let $\Phi : SL_3(\Bbb C) \to  G_2(\Bbb C)$ be
the standard inclusion of the L-groups of $PGL_3(F)$ and 
$G_2(F)$. Let $Y$ be a tempered spherical representation of 
$PGL_3(F)$. Let $s\in SL_3(\Bbb C)$ be its Satake parameter. Let
$\tilde Y$ be the tempered spherical representation of $G_2(F)$ 
whose Satake parameter is $\Phi(s)$. Then 
$$
\Theta (Y) = \{\tilde Y\}. 
$$
\endproclaim 

\demo{Proof} Write $Y=\Ind_{Q}^{GL_3}X$ (note that there are up to 
three different choices for $X$). Then the restriction
of $X$ to $GL_2(F)$ (a factor of $L$) is the irreducible spherical 
representation whose parameter is 
$(\chi |\cdot |^{1/2}, \mu |\cdot |^{1/2})$ where $\chi$ and $\mu$ are unitary
characters. Let  
$\tilde {X}$ be the irreducible spherical representation of $GL_2(F)$ 
whose parameter is  
$(\chi^{-1} |\cdot |^{3/2}, \mu^{-1} |\cdot |^{3/2})$. 
Note that $\ca(A)$ is the regular representation of
$GL(2)$. After taking into account the twist with $|\det |^2$ it 
follows that there is a unique quotient of $\ca(AA)$ 
isomorphic to $\tilde X\otimes X$ as a $GL_2(F)\times GL_2(F)$-module. Here 
the first copy of $GL_2(F)$ is the Levi factor of $P_2$ and the second is
the factor of $L$. Therefore,   
$$
\Ind_{P_2}^{G_2}(\tilde X)\otimes X
$$
is a quotient of $\ca(AA)$. Write  
$ \tilde Y =\Ind_{P_2}^{G_2}\tilde X$. Let $U_2$ be the unipotent 
radical of $P_2$. It is Heisenberg group. Let $Z$ be the center of $U_2$. 
The action of $GL_2(F)$ (Levi factor of $P_2$) on $Z$ is $\det$ and  
on $U_2 /Z$ it is isomorphic to $S^3(F^2)\otimes \det^{-1}$. It 
follows that $\rho_{U_2}=|\det|^3$. Hence $\tilde Y$ is tempered.  
It is irreducible by a result of Keys 
\cite{Ke}.  
We have to show that $\tilde Y\otimes X$
is a quotient of $V_U$ and that $X$ does not appear as a quotient 
composition factors of $V_U$ different from $\ca(AA)$. 

We again use Lemma 2.5, but in this 
case $T$ will be in the Bernstein center (see \cite{BD})
of $GL_2(F)\subseteq L$.
Recall that the component of the Bernstein center corresponding 
to representations generated by their Iwahori-fixed vectors is isomorphic 
to 
$$
\Bbb C[x,x^{-1},y,y^{-1}]^{W}
$$
where $W=\{1,w\}$, $w(x)=y$ and $w(y)=x$ is the Weyl group of $GL_2$. 
Let $I$ be the Iwahori subgroup of $GL_2(F)$. 
Let $\varpi$ be the uniformizing element in $F$. Then any unramified 
character $\chi$ of $F^*$ is determined by its value on $\varpi$. 
If $E$ is a subquotient of an induced representation with the parameter
$(\chi_1, \chi_2)$ then 
$$
(x+y) = \chi_1(\varpi) + \chi_2(\varpi) \text { and } 
xy = \chi_1(\varpi)\chi_2(\varpi) 
$$
on $E$. 
Let
$$
T_1=q^{3/2}(x^{-1}+y^{-1})-(xy)^{-1},
$$
where $q=|\varpi|^{-1}$. 
On $X$  $T_1$ acts as the scalar  
$$
q^2(\chi(\varpi)^{-1} +\mu(\varpi)^{-1})-q\chi(\varpi)^{-1}\mu(\varpi)^{-1}, 
$$
which if real, is strictly less then $2q^2$. On the other hand, 
$\ca(B)$ is the regular representation of $GL_1$, so
$\ca(B)$ as a $GL_2(F)$-module consists of induced representations
whose inducing parameters are $(|\cdot|^{3/2}, \chi)$. 
It has an Iwahori-fixed vector only when $\chi$ is an unramified 
character. On such induced representations $T_1$ acts as 
$$
q^3=q^{3/2}(q^{3/2}+\chi(\varpi)^{-1})-q^{3/2}\chi(\varpi)^{-1}. 
$$
Hence $T_1$ acts on $\ca(B)^{I}$ as the scalar $q^3$, which is always 
bigger then $2q^2$. Therefore $X$ can not be a quotient of $\ca(B)$.  

Let 
$$
T_2 =xy.
$$
Then $T_2$ acts on $X$ as  
$|\varpi|\chi(\varpi)\mu(\varpi)$ which is  
different from  $|\varpi|^{2}$ and $|\varpi|^4$, the eigenvalues of $T_2$ on 
$V^I_N$. Therefore $X$ can not be a quotient of $V_N$. Finally, Lemma 2.5 
applied to $T_1$ and $T_2$ shows that 
$\tilde Y\otimes X$ is a quotient of $V_U$. 

To finish the proof we have to check that the correspondence $Y\to 
\tilde Y$ is Langlands. Since $\tilde Y =\Ind_{P_2}^{G_2}\tilde X$ 
and $P_2$ is spanned by a short root, the Satake parameter of 
$\tilde Y$ sits in a Levi factor of a parabolic subgroup of 
$G_2(\Bbb C)$ spanned by a long root. Since $SL_3(\Bbb C)$ is 
spanned by long roots of $G_2$, $\tilde Y$ must be a lift of 
a representation of $PGL_3(F)$ induced from $Q$: $Y$ or $Y^*$. 
Note that replacing the pair $(\chi,\mu)$ by 
$(\chi^{-1},\mu^{-1})$ does not change $\tilde Y$
but replaces $Y$ by $Y^*$. The theorem is 
proved. 
\enddemo

\head 5. Dual pair $G_2(F) \times PGSp_6(F)$ \endhead 

Let $G$ be the simply connected group of type  $E_7$. Let $P=MN$ 
be the maximal parabolic subgroup of $G$ such that $M$ is the 
reductive group of type $E_6$ discussed in the last section. 
The nilpotent radical $N$ is commutative and isomorphic to $J$ 
as an $M$-module. Let  
$G_2(F) \times GL_3(F)$ be the dual pair in $M$, described in Section 4. 
The centralizer of $G_2(F)$ in $G$ is $Sp(6)$. This can be easily seen on 
the level of Lie algebras. Let 
$\frak g$ be the Lie algebra of $G$. 
Then 
$$
\frak g = \bar \frak n + \frak m + \frak n
$$  
and write 
$$
C_{\frak g}(G_2(F))= 
 \bar \frak u + \frak l + \frak u 
$$
where $\frak l \subset \frak m$, $\overline{\frak u} \subset \overline{\frak n}$ and
$\frak u \subset \frak n$. Then $\frak l = gl(3)$ and $\frak u 
\subset \frak n$ corresponds to the inclusion $J_F \subset J$ 
since 
$$ J^{G_2(F)} = J_F.$$ 
Therefore  $C_{\frak g}(G_2(F))= sp(6)$ whose Siegel parabolic 
subalgebra is $\frak l + \frak u$. Note also that we have a 
distinguished $SL_2$, 
$$ SL_2(F) 
\subset  Sp_6(F) \subset  G
$$
such that if $(e,h,f)$ is the standard basis of
$sl(2)$ (the Lie algebra of $SL_2$) then 
$h \in \frak m$, $f \in \overline{\frak n}$
and
$ e\in \frak n$. Moreover $e$ corresponds to the unit 
element of $J$ under the identification of $\frak n$ and $J$. 
Let $G_A$ be the split adjoint group of type $E_7$, and let $P_A=M_A N$ 
be the corresponding maximal parabolic subgroup. Then $M_A$ can be 
defined as the group of isogenies of the Dickson form. 
We have a dual pair 
$G_2(F) \times PGSp_6(F) \subset G_A$. 
Note that $PGSp_6(F)$ is 
generated by the distinguished $PGL_2(F)\subset  PGSp_6(F)$ and the image of
$Sp_6(F)$. 

Let $Q=GL_3(F)U$ be the Siegel parabolic subgroup of 
$Sp_6(F)$ and $Q_A=LU$ be the corresponding maximal parabolic of $PGSp_6(F)$.  
Let $V$ be the minimal representation of $G_A$. We want
to compute $V_U$. Since
$$
0 \to {\Cal C}_{c}^{\infty}(\oz)_U \to V_U  
\to V_N \to 0 
$$ 
we need to understand $\ca(\oz)_U$. Let $NN$ and $\overline{NN}$ be the 
orthogonal complements of $\overline{U}$ and $U$ in $N$ and 
$\overline{N}$. 
Since $U$   
can be identified with $J_F$,  
it follows that 
$$
NN = \{ 
\pmatrix 
0 & x & y \\
\overline{x} & 0 & z \\
\overline{y} & \overline{z} & 0 
\endpmatrix 
|~ \overline{x}=-x, \overline{y} =-y \text { and } \overline{z} =-z \}.
$$  
Let $\om =\oz \cap \overline{NN}$. As in Lemma 2.2  
$\ca(\oz)_U = \ca(\om )$ and we have:

\proclaim{Proposition 5.1} Identify $\overline{NN}$ with the set of triples of
traceless Octonions $(x,y,z)$. Let $f\in \ca(\om)$. The action of $G_2(F)\times L$ is
given by 
\roster
\item
$$
\pi(g)f((x,y,z)) = f((g^{-1}x, g^{-1}y, g^{-1}z)), \qquad g\in G_2(F). 
$$
\item
$$
\pi\pmatrix a & \\ & b \\ \endpmatrix f((x,y,z)) = |\frac{a}{b}|^6
f(\frac{a}{b}(x,y,z)), \qquad
\pmatrix a & \\ & b \\ \endpmatrix\in PGL_2(F)\subset PGSp_6(F).
$$
\item
$$
\pi(g)f((x,y,z)) = |\det g|^4
f(\det g(x,y,z)\,^tg^{-1}), \qquad g \in GL_3(F).
$$
\endroster
\endproclaim

\demo{Proof} It follows from Theorem 1.1. Note, however, that in this case the formulas 
are describing the action on $\ca(\om)$. 
\enddemo 

To keep our notation simple we work in $G$. 
We have to understand the structure of $G_2(F) \times GL_3(F)$
orbits on $\om$ which is, of course, the
same as the structure of orbits on $\omega\omega=\omega \cap NN$.

\proclaim{Proposition 5.2} 
\roster 
\item $\om = \{(x,y,z)|\overline{x}=-x,  
 \overline{y} =-y,  \overline{z} =-z \text { and }x^2= y^2=z^2=xy=
 xz=yz=0 \}$.
\item Let $AA$ and $BB$ be the subsets  
of $\om$ consisting of 
all triples $(x,y,z)$ such that the space $Fx+Fy+Fz$ has dimension $2$ and $1$ 
respectively. 
Then $AA$ and $BB$ are $G_2(F)\times GL_3(F)$-orbits and 
$$\om =  
 AA \cup BB. $$

\endroster 
\endproclaim 

\demo{Proof} Let us prove the corresponding statement for 
$\omega\omega$. Let   
$$
n= 
\pmatrix 
0 & x & y \\
-x & 0 & z \\
-y & -z & 0 
\endpmatrix \in \omega\omega .
$$
Since $n$ is a traceless matrix in $J$ it will be singular iff 
$n^2=0$. But this is equivalent to $x^2=y^2= 
z^2=xy=xz=yz = 0$. The first part of the proposition is proved.

Again, recall that $G_2$ has three orbits on the set of 
spaces 
of traceless Octonions with the property that the Octonion
multiplication is trivial. These are characterized by their dimension;
the possible choices being $0,1,2$.
The stabilizers of the nontrivial spaces are the maximal parabolics
of $G_2$. It follows that $x$, $y$ and $z$ are linearly dependent. 
hence $\om = AA \cup BB$. It remains to show that $AA$ and $BB$ are 
single orbits. The proof is analogous to the proof of Proposition
4.2. 

Let $x$ be a traceless Octonion such that $x^2=0$. Let
$P_1$ be the maximal parabolic subgroup of $G_2(F)$ stabilizing
the line $Fx$. Consider
$$
B=\{ (ax,0,0)|a\in F \text { and } a\neq 0\}. 
$$
Let $Q_1$ be the maximal parabolic of $GL_3(F)$ stabilizing $B$. 
Then $P_1 \times Q_1$ acts transitively on $B$ and 
$$
BB=  (G_2(F)\times GL_3(F)) \times_{(P_1\times Q_1)} B.
$$
Let $y$ and $z$ be two traceless and linearly 
independent Octonions such that $y^2=z^2=yz=0$. Let  
$P_2$ be the maximal parabolic subgroup of $G_2(F)$ stabilizing
the space $Fy+Fz$. Consider
$$
A=\{ (ay+bz,cy+dz,0)|a,b,c,d\in F {\text { and }} ad-bc\neq 0\}. 
$$
Let $Q_2$ be the maximal parabolic subgroup of $GL_3(F)$ stabilizing $A$. 
Then $P_2 \times Q_2$ acts transitively on $A$ and 
$$
AA =  (G_2(F)\times GL_3(F)) \times_{(P_2 \times Q_2)} A.
$$
The proposition is proved.  
\enddemo 

We can now summarize the structure of $V_U$ as a $G_2(F) \times 
GL_3(F)$-module. 

\proclaim{Theorem 5.3}
$V_U$ has a filtration with three subquotients:
$$
\ca(AA),\,\, \ca(BB),\, \text { and } V_N.
$$
Moreover 
\roster 
\item 
$$
\ca(AA) = \ind_{P_2\times Q_2}^{G_2\times GL_3}(\ca(A))\otimes |\det|^4  
$$ 
\item  
$$
\ca(BB)=\ind_{P_1\times Q_1}^{G_2\times GL_3}(\ca(B))\otimes |\det|^4
$$
\item 
$$
V_N \cong V_M \otimes |\det |^2 + 1\otimes |\det|^4 
$$
\endroster 
as $G_2(F) \times GL_3(F)$-modules. Here $\det$ denotes the usual 
determinant of $3\times 3$ matrices. 
\endproclaim
 
We are now ready to state and prove the main theorem. 

\proclaim{Theorem 5.4} Let $\Phi :  G_2(\Bbb C) \to Spin_7(\Bbb C)$ be
the standard inclusion of the L-groups of $G_2(F)$ and 
$PGSp_6(F)$; $G_2(\Bbb C)$ fixes a non-zero vector in the $8$-dimensional 
spin representation of $Spin_7(\Bbb C)$. 
Let $Y'$ be a tempered spherical representation of 
$PGSp_6(F)$. 
Then $\Theta(Y')$ is not empty only if the Satake parameter of 
$Y'$ is $\Phi(s)$ for some $s$, a Satake parameter of a 
tempered spherical representation $\tilde Y$ of $G_2(F)$. In that 
case 
$$
\Theta (Y') = \{\tilde Y\}. 
$$
\endproclaim 

\demo{Proof} 
Let $X'$ be a spherical tempered representation of $PGSp_6(F)$. 
Recall that every tempered spherical representation of $PGSp_6(F)$ is 
fully induced (see \cite{T1}). 
To simplify notation we work with $Sp_6(F)$. 
The restriction of $X'$ to $Sp_6(F)$ can be written as  
$\Ind _{Q}^{Sp_6}X\otimes |\det|^2$ where $X$ is a tempered spherical 
representation of $GL_3(F)$ (note that $\rho_U=|\det|^2$).

Let $Y'$ be a spherical tempered representation of
$PGSp_6(F)$ whose Satake parameter is $\Phi(s)$. 
This means that the restriction of $Y'$ to $Sp_6(F)$
can be written as  $\Ind_{Q}^{Sp_6}(Y\otimes |\det|^2)$
where $Y$ is a tempered spherical representation of $PGL_3(F)$. 
Also  
$\tilde Y$, 
the representation of $G_2(F)$ whose parameter is $s$, is the lift 
of $Y$ from $PGL_3(F)$. 
By Theorem 5.3, $V_M\otimes |\det|^2$ is a quotient of $V_U$ ($M$ 
is a group of type $E_6$). 
It follows from the Frobenius reciprocity that $\tilde Y\in \Theta(Y')$. 

The rest of the theorem follows form the knowledge 
of $V_U$. Indeed, let $X'$ be a tempered spherical representation 
of $PGSp_6(F)$. Again, it is a fully induced representation  
so its restriction to $Sp_6(F)$ can be written as  
$\Ind _{Q}^{Sp_6}X\otimes |\det|^2$ where $X$ is a tempered spherical 
representation of $GL_3(F)$. 
If $\tilde Y\otimes X'$ is 
a quotient of $V$ for some irreducible representation $\tilde Y$ of  
$G_2(F)$, then $\tilde Y \otimes (X\otimes |\det|^2)$ 
is a quotient of $V_U$, i.e.  
it is a quotient of one of the three pieces in Theorem 5.3. 
For example, if it is a quotient of $\ca(AA)$ then  
$$
\tilde Y\otimes (X\otimes |\det|^2) =
\Ind_{P_2}^{G_2}(\tilde X)\otimes(X\otimes |\det|^2) 
$$ 
where the parameter of $\tilde X$ 
is $(|\cdot|^{3/2}\chi^{-1}, |\cdot|^{3/2}\mu^{-1})$ and 
the parameter of $X$ is $(\chi,\mu,
\chi\mu)$ for some unitary characters $\chi$ and $\mu$.
Since  $\Ind X \otimes |\det|^2  
= \Ind Y\otimes |\det|^2$, where $Y$ is a representation of 
$PGL_3(F)$ with parameter $(\chi,\mu,\chi^{-1}\mu^{-1})$, 
we are back in the 
situation as in the beginning of the proof. The theorem is 
proved. 

\enddemo 

With only minor modifications of the preceeding proof one 
can easily show the following proposition. 

\proclaim{Proposition 5.5} Let $Y$ be a spherical representation 
of $PGSp_6(F)$ which appears as a quotient of $V$, the minimal 
representation of $E_7$. Then the Satake parameter of $Y$ is 
$\Phi(s)$ for some $s\in G_2(\Bbb C)$. 
\endproclaim 

\head 6. Some facts about Jordan algebras  \endhead

In this section we shall describe structure of exceptional groups 
using simple Jordan algebras of rank $3$.  

Let $D$ be a non-degenerate subalgebra of Octonions. 
Let $J_D \subseteq J$ be the corresponding Jordan algebra of rank $3$. 
The Dickson form on $J_D$ is given by 
$$6(x,y,z)= 
2tr(xyz) - tr(x)tr(yz)-tr(y)tr(zx)-tr(z)tr(xy)-tr(x)tr(y)tr(z).
$$
Define 
$$
\frak l_D = \{ x\in gl(J_D) \,|\, (xv,v,v)=0 \text { for all } v\in J_D\}.
$$

Now assume that $D=F$, $F+F$, $M_2(F)$ ($2\times 2$-matrices) or $\Bbb O$.
Let $G_D$ be a simply connected group of type $F_4$ or $E_n$, $n=6,7$ or $8$ 
respectively. 
Let $\frak g_D$ be its Lie algebra. Let $\Delta$ be the Dynkin diagram
of $\frak g_D$. We shall identify it with a set of simple roots.  
Let $\tilde \alpha$
be the highest positive root. 
Let $\alpha$ be the unique simple root not perpendicular to 
$\tilde \alpha$. Let $\frak m_D\subset \frak g_D$ be the simple 
algebra whose Dynkin diagram is $\Delta \setminus \{\alpha\}$.
Extend $\Delta$ 
by adding $-\tilde \alpha$. 
Let $\beta$ be the unique simple root not perpendicular 
to $\alpha$. Remove the vertex corresponding 
to the simple root $\beta$. The extended diagram  
breaks into several pieces, one of which is an $A_2$ diagram 
corresponding to $\{\alpha, -\tilde \alpha\}$. 
Then $\frak l_D\subset \frak g_D$ corresponds to the 
rest of the diagram. Consider the adjoint action of  
${sl}(3) +\frak l_D$ on $\frak g_D$. 
One has a decomposition
$$
\frak g_D=sl(3)+\frak l_D + V \otimes J_D + V^* \otimes J^*_D 
$$ 
where $V$ is the standard $3$-dimensional representation of 
$sl(3)$ (see \cite{HP}).
Let $M_D$ and $L_D$ be the corresponding simply connected 
groups and let  $H_D$ be
a subgroup of $L_D$ centralizing the identity of $J_D$. The
possible choices are
given by the Freudenthal's magic square:
$$
\matrix
G_D & M_D & L_D & H_D & J_D\\
F_4 & Sp_6 & SL_3 & SO_3 & S^2 V \\ 
E_6 & SL_6 & SL_3 \times SL_3 & SL_3 & V\otimes V^* \\
E_7 & E_6 & SL_6 & Sp_6 & \wedge^2 F^6 \\
E_8 & E_7 & E_6 & F_4  &  F^{27}
\endmatrix
$$
In the above table we have identified $J_D$ with a "usual"   
representation of $L_D$ whenever possible. Since the 
centralizer of $H_D$ in $J_D$ is one-dimensional, it follows 
that the centralizer of $H_D$ in $\frak g_D$ is (see \cite{FH}) 
$$
\frak g_2 = sl(3) + V + V^*.
$$

Assume that $D\neq F$.
We shall need few facts about the action of $H_D$ on $J_D^0$, 
the set of traceless elements in $J_D$.  
Recall that a point $x \in J^0_D$ is  
singular if $x^2=0$.  
The group $H_D$ acts transitively 
on the set of singular vectors in $J_D^0$. Denote by $Q_1$  
a parabolic subgroup of $H_D$ stabilizing   
a singular line in $J_D^0$.  

Let $Fy +Fz \subset J^0_D$ be a singular space, i.e. 
each point is singular. This is equivalent to 
$$
y^2=z^2=yz=0.  
$$
We need to know $H_D$-orbits of  
singular two-dimensional subspaces of 
$J_D^0$. We have two different cases. 

If $\dim D=2$ then 
$J_D$ is isomorphic to the algebra of all $3\times 3$-matrices
with coefficients in $F$ and $H_D=SL_3$ acts by conjugation. The
restriction of the Dickson form to $J_D$ is given by
$$
(x,x,x)=\det(x)
$$
where $\det$ is the determinant of $3\times 3$ matrices.
In this case, a singular point in $J_D^0$ is a nilpotent 
rank-one matrix. There are two $SL_3$-orbits of singular 
two-dimensional spaces. Let $Fy+Fz$ be a singular space. Then 
either the images or the kernels of $y$ and $z$ coincide. 
Denote by $Q^+$ or $Q^-$ the stabilizer of the singular space 
depending on the orbit. It is a maximal parabolic subgroup of 
$SL_3$. 

In other cases, a stabilizer of a singular two-dimensional 
space is a parabolic subgroup only if the space is "amber"
(\cite{A}).

\proclaim{Definition 6.1}   
Let $S\subset J_D^0$ be a singular two-dimensional space. We say that $S$ is 
amber if $S\subset xJ_D$ for every $x\in S$. 
\endproclaim  

\proclaim{Proposition 6.2}
If $\dim D \geq 4$ then $H_D$ acts transitively on the 
space of amber two-dimensional subspaces of $J_D^0$.  
\endproclaim 

\demo{Proof}
If $\dim D=8$ this is a result of Aschbacher, 9.3-5 \cite{A1}. 
Let $Fy +Fz=S$ be an amber space. Since $H_D$ acts transitively 
on the set of singular vectors, we can assume that $z=x$, where $x$ 
is a fixed singular vector in $J^0_D$. Let $Q_1$ be a parabolic 
subgroup of $H_D$ stabilizing the line $Fx$. Let $L_1$ be a 
Levi factor of $Q_1$. Let $I_x \subset J_D$ be an $L_1$-invariant 
space such that 
$$
xJ_D \cap J_D^0 = Fx + I_x .
$$
To finish the proof, we have to show that $L_1$ acts transitively on 
the set of singular points in $I_x$.

If $\dim  D=4$ then $L_1= GL_2 \times SL_2$ and $I_x$ can be identified with the 
space of $2\times 2$-matrices. 
If $\dim D=8$ then 
$[L_1,L_1]=Spin_7$ and $I_x$ is the corresponding $8$-dimensional spin-representation. 
In both cases $L_1$ has two non-trivial orbits in $I_x$, the smaller consisting of singular
points. The proposition is proved.  
\enddemo 

Let $Q_2$ be the stabilizer of an amber plane.  
We finish this section with a description of $Q_2$ in both cases.  
Identify the amber space with $F^2$. Let $L_2$ be a Levi factor of $Q_2$.  

If $\dim D=4$ then $Q_2$ is contained in a Siegel parabolic subgroup of $Sp(6)$ and 
$L_2=GL_2 \times GL_1$. 
The action of $L_2$ on the 
amber space is the standard action of $GL_2\times GL_1$ on the space of $2\times 1$
matrices.

If $\dim D=8$. Then $L_2$ is a quotient of $GL_3\times GL_2$ obtained by 
identifying centers of both factors (consider  
$GL_3\times GL_2$ acting 
on $3\times 2$ matrices, for example). Let $(l_3 , l_2) 
\in GL_3\times GL_2$. The action of $L_2$ on the 
amber space is given by 
$$
\det(l_3)\det(l_2) \l_2 x\,\, x\in F^2.  
$$
In particular, the action of $L_2$ on the amber space induces a
non-split exact sequence 
$$
1\to SL_3 \to L_2 \to GL_2 \to 1. 
$$ 
We have $[L_2,L_2]=SL_3 \times  SL_2$ where $SL_3$ is spanned by 
long roots and $SL_2$ by short roots.

\head 7. Jacquet functor for $G_2$ - Heisenberg parabolic \endhead

We continue to use the notation and definitions introduced in the 
previous section. 
Let 
$$
x= 
\pmatrix 
1 & 0 & 0 \\
0 & 0 & 0 \\
0 & 0 & -1 
\endpmatrix 
\in sl(3)\subset \frak g_D. 
$$
Let $\frak g(k) = \{ y\in \frak g_D | [x,y]=ky\}$. Then $\frak g(k) 
\neq 0$ for $k=-2,-1,0,1,2$. Let $\frak p = \frak m + \frak n$ 
be the maximal parabolic subalgebra with $\frak m=\frak g(0)$ and 
$\frak n = \frak g(1) + \frak g(2)$. Note that $[\frak m, \frak m] 
= \frak m_D$ and $\frak n$ is 
a Heisenberg Lie algebra with center $\frak g(2)$.
Let $P=M N$ be the corresponding maximal parabolic subgroup. 
Let $Z$ be the center of $N$ and let $N_D$ be the quotient of
$N$ by $Z$. Then $N_D$ is an irreducible $M_D$-module. 
Since $N_D\cong \frak g(1)$, the restriction to $L_D$ is isomorphic to
$$
F + J_D + J_D^* + F.
$$

The centralizer of $H_D$ in $G_D$ is $G_2(F)$, and
$$
G_2(F) \cap P = P_2 = GL_2(F) U_2
$$
is the Heisenberg parabolic of $G_2(F)$.
Then $Z \subset U_2$ and let $U$ be the
quotient of $U_2$ by $Z$.

Let $\Omega$ be the set of singular vectors in $N_D$.
Let $\overline{N}$ denote the opposite
nilpotent subgroup. Let $\overline{N}_D$, $\overline{Z}$...
be the corresponding objects attached to $\overline{N}$.
Note that the Killing form on $\frak g_D$, the Lie algebra of
$G_D$, defines a non-degenerate pairing $<\cdot ,\cdot>$ 
between $N_D$ and
$\overline{N}_D$.

\proclaim{Theorem 7.1} Assume that the residual characteristic of
$F$ is odd. Let $V$ be the minimal representation of 
$G_D$. Let $Z$ be the center of $N$ as above. Then $V_Z$ 
has a $P$-invariant filtration
$$
0 \to {\Cal C}_{c}^{\infty}(\Oz) \to V_Z  
\to V_N \to 0
$$
where $\Cal C_{c}^{\infty}(\Oz)$ denotes the space of locally
constant, compactly supported functions on $\Oz$. 
\roster
\item Let $f\in \Cal C_{c}^{\infty}(\Oz)$. Then 
$$\
\pi(n)f(x) = \psi(<x,n>)f(x), \qquad n\in N
$$
and
$$
\pi(m)f(x) = |\det(m)|^{\frac{s}{d}}f(m^{-1}x), \qquad m\in M. 
$$
\item
$$
V_N \cong V_M \otimes |\det|^{\frac{t}{d}} + |\det|^{\frac{s}{d}}  
$$
where $V_M$ is the minimal representation of $M$ (center acting
trivially).
\endroster
In the above formulas $\det$ is the determinant of
the representation of $M$ on $N_D$, $d$ is the dimension of $N_D$ and
the values of $s$ and $t$ are given in the following table.
$$
\matrix
G & s & t \\
E_6 & 4 & 3 \\
E_7 & 6 & 4 \\ 
E_8 & 10 & 6 
\endmatrix
$$
\endproclaim

\demo{Proof} Part (2) is Prop. 4.1 \cite{S}.  
Let $W$ be the kernel of the projection of $V_Z$ onto $V_{N}$. 
Let $x\in \overline N_D$ and let 
$$
\psi_x(y)=\psi(<x,y>) 
$$
be a character of $N$. The theorem of Moeglin and Waldspurger \cite{MW} 
implies that $\dim W_{N,\psi_x}=0$ or $1$ and it is one 
if and only if $x\in \Oz$. Let $x\in \Oz$. Let $M_x$ be the 
stabilizer of $x$ in $M$ and $\chi$ the character of $M_x$ 
describing the action of
$M_x$ on $W_{N,\psi_x}$. By Frobenius reciprocity there 
exists a non-trivial $P$-homomorphism 
$$
T : W \to  \Ind_{M_x N}^{P}(\chi \otimes \psi_x). 
$$
Let $\Cal C^{\infty}(\Oz)$ denote the space of locally
constant functions on $\Oz$. 
Note that we have an inclusion 
$$
\Ind_{M_x N}^{P}(\chi \otimes \psi_x) \subseteq 
\Cal C^{\infty}(\Oz). 
$$
Let $w\in W$ and $f=T(w)$. We need to show that $f$ 
is a compactly supported function on $\Oz$. 
Let $N_{D,k}$, $k \in \Bbb Z$ be a family of  lattices of $N_D$
such that $\cup_k N_{D,k}= N_D$ and $\cap_k N_{D,k}=0$.
Let $\overline N_{D,k}$ be their dual lattices in $\overline N_D$.
Since $W$ is a smooth module, there exists a small integer $k$
depending on $w$ such that $\pi(n)w=w$ for all $n\in N_{D,k}$. 
This implies that $f$ is supported inside $\overline N_{D,k}$. 
Since $W_{N}=0$ there exists an integer $l$ depending on $w$
such that 
$$
\int_{N_{D,l}}\pi(n)w dn = 0
$$
(see 2.33 \cite{BZ}). 
This implies that $f$ is supported outside $\overline N_{D,k}$.
Since $\Oz$ is locally closed and the boundary is $\{ 0\}$, it 
follows that $f\in \Cal C_{c}^{\infty}(\Oz)$. Let $\ind$ denote 
smooth induction with compact support. By the Bernstein-Zelevinsky
analogue of Mackey Theory (see \cite{BZ}, pages 46-47)  
$$
\ind_{M_x N}^{P}(\chi \otimes \psi_x)= \Cal C_{c}^{\infty}(\Oz)
$$ 
is irreducible $P$-module. Hence $T(W)=\Cal C_{c}^{\infty}(\Oz)$. 
Let $W'$ be the kernel of $T$. Since 
$$
\dim W_{N, \psi_x}=\dim \Cal C_{c}^{\infty}(\Oz)_{N, \psi_x}
$$
for any $x\in \overline N_D$, it follows that $W'_{N, \psi_x}=0$ 
for any $x\in \overline N_D$ (2.35 \cite{BZ}). Therefore 
$W'=0$ by 5.14 \cite{BZ}. 

Note that the inclusion $M_x \to M$
induces an isomorphism $M_x/[M_x ,M_x] \cong M/[M,M]$ This can be 
easily checked by choosing $x$ to be in $\frak g_{D,\alpha} \subset 
\frak g(1)$, $\alpha$ is defined in the previous section. Hence  
$\chi$ is a character of $M$ and to finish the proof we have to 
show that 
$$
\chi(m) =|\det(m)|^{\frac{s}{d}} \qquad m\in M. 
$$
Let $GL_2(F)$ be the Levi factor of $P_2$ defined above. The 
inclusion of $GL_2(F)$ into $M$ induces in isomorphism 
$GL_2(F)/SL_2(F) \cong M/[M,M]$. This can be seen by looking 
at the action of $GL_2(F)$ and $M$ on $Z$. Therefore it 
suffices to find the restriction of $\chi$ to $GL_2(F)$.
In Section 8 we shall use the information on correspondences obtained in 
previous sections to find the character. 
\enddemo

Let $NN$ be
the orthogonal complement of $\overline{U}$ in $N_D$. Let $J_D^0$ 
be the set of traceless elements of $J_D$.  
Then $NN$ as a $GL_2(F)\times H_D$-module is isomorphic to 
$$
 F^2 \otimes J_D^0. 
$$ 
Write 
$$
\Omega\Omega = \Omega \cap NN. 
$$ 

\proclaim{Proposition 7.2} 
\roster 
\item $\Om  = \{(y,z)\neq (0,0)|~ y,z  
\in J_D^0 \text { and the space } Fy+Fz \text { is amber } \}$.
\item Let $AA$ and $BB$   
be the subsets of $\Om$ consisting 
of all pairs $(y,z)$ such that the space $Fy+Fz$ has dimension $2$ and $1$  
respectively. Then $BB$ is a $GL_2(F)\times H_D$-orbit. If $\dim D\geq 4$
then $AA$ is a $GL_2(F)\times H_D$-orbit. If $\dim D=2$ then $AA$ is a 
union of two orbits.
\endroster 
\endproclaim 

\demo{Proof} 
Recall that a two-dimensional space $S \subset J_D^0$ is  
amber if every element is singular and 
$$
S\subset xJ_D  
$$
for all $x\in S$.
As before, write    
$$
N_D = F + J_D +J_D^* +F. 
$$
Let $Q_D=LU_D \subseteq M_D$ be the stabilizer of the last summand in the above 
decomposition. Then $[L,L]=L_D$, $U_D\cong J_D$ as an $L_D$-module
and $Q_D$ stabilizes the partial flag 
$$ 
N_D \supset J_D +J_D^* + F \supset J_D^* +F \supset F.
$$ 
More precisely, let $u\in U_D$, $(a,y,z,b)\in N_D$ and 
$(a',y',z',b')=u(a,y,z,b)$. Then 
$$
a'=a,\, y'= y + au,\, z'=z+2u\times y + au\times u,
$$
$$
b'=b+tr(uy)+3(u,u,z)+a\det u,  
$$ 
where $\det u=(u,u,u)$ and $u\times y$, the cross product, is 
an element of $J_D$ such that 
$$
tr((u\times y)x) = 3(x,y,u)
$$
for all $x\in J_D$ (see \cite{Ki} page 143). 

\proclaim{Lemma 7.3} The group $Q_D$ has $4$ orbits on $\Omega$.
Their representatives are 
$$
v_1=(1,0,0,0),\, v_2=(0,x,0,0)
$$
$$
v_3=(0,0,x,0),\, v_4=(0,0,0,1)
$$
where $x$ is any singular element in $J_D$. 
\endproclaim 

\demo{Proof} Note that $\Bbb P(\Omega)=M_D/Q_D$. We have to compute 
$Q_D\backslash M_D/Q_D$ which is the same as $W_L\backslash W_M/W_L$,
here $W_M$ and $W_L$ denote the Weyl groups of $M_D$ and $L_D$. 
Since $N_D$ is a miniscule representation i.e. weight vectors are 
all contained in one $W_M$-orbit, it follows that they are 
parametrized by $W_M/W_L$. On the other hand $J_D$ is a miniscule 
representation of $L_D$ so $W_L\backslash W_M/W_L$ has four orbits 
and it follows easily that $v_i$, $(1\leq i \leq 4)$ form 
a complete set of representatives of $Q_D$-orbits. The lemma
is proved.
\enddemo 

Let $(0,y,z,0) \in \Om$. If $y\neq 0$ then Lemma 7.3 implies that it is $Q_D$ 
conjugated to $v_2$. Therefore $y$ is in the $L$-orbit of $x$, hence $y$ 
is singular. Since the action of $GL_2(F)$ is 
$$
(y,z)
\pmatrix 
a & b \\
c & d 
\endpmatrix 
=(ay+cz,by+dz) 
$$
the same argument implies that any
element of $Fy+Fz$ is singular. 

Note that $(0,y,z,0)$ is $Q_D$, hence $U_D$-conjugated to  
$(0,y,0,0)$. But this 
means that $z=2 u \times y$ for some 
$u\in U_D\cong J_D$ such that $tr(uy)=0$. 
Since  
$
tr((u\times y)v)=3(u,y,v)  
$
for all elements $v\in J_D$ and 
$$
6(u,y,v)=2tr(uyv)-tr(u)tr(yv)
$$
(use $tr(y)=tr(uy)=0$), it follows that 
$$
tr((2u\times y)v)=tr(u^*yv)
$$
for all $v\in J_D$, here $u^*=2u-tr(u)$. Therefore 
$z=u^*y\in yJ_D$. Since the same argument can be repeated for any
linear combination of $y$ and $z$, the first part of the 
proposition follows. 

Let $x \in J_D^0$ such that $x^2=0$. Let
$Q_1$ be the parabolic subgroup of $H_D$ stabilizing
the line $Fx$. Consider
$$
B=\{ (ax,0)|a\in F \text { and } a\neq 0\}. 
$$
Let $QQ\subset GL_2(F)$ be the Borel subgroup stabilizing 
the line $B$. 
Then $QQ\times Q_1$ acts transitively on $B$ and 
$$
BB=  (GL_2(F)\times H_D) \times_{(QQ \times Q_1)} B.
$$
Assume now that $\dim D\geq 4$. Let $y,z\in J_D^0$ such that the space $Fy+Fz$ 
is amber.
Let $Q_2$ be the parabolic subgroup of $H_D$ stabilizing
the space $Fy+Fz$. Consider
$$
A=\{ (ay+bz,cy+zd)|a,b,c,d\in F \text { and } ad-bc\neq 0\}. 
$$ 
Then $GL_2(F)\times Q_2$ acts transitively on $A$ and 
$$
AA =  H_D \times_{Q_2} A.
$$
If $\dim D=2$ then we have two orbits of singular two-dimensional 
spaces. Let $Fy^++Fz^+$ and $Fy^-+Fz^-$ be their representatives and $Q^+$ 
and $Q^-$ their stabilizers in $SL_3(F)$. One can define $A^+$ and 
$A^-$ as above, hence    
$$
AA =  SL_3(F) \times_{Q^+} A^+ \cup SL_3(F) \times_{Q^-} A^-.
$$
The proposition is proved. 
\enddemo 
Since $\ca(\Oz)_{U_2}=\ca(\Om)$, the structure of 
$V_{U_2}$ as $GL_2(F) \times 
H_D$ is given as follows. 

\proclaim{Theorem 7.4}
Assume that the residual characteristic of $F$ is odd. 
$V_{U_2}$ has a filtration with three quotients:
$$
\ca(AA),\,\, \ca(BB),\, \text { and } V_{N}.
$$
Moreover,  
as $GL_2(F) \times H_D$-modules,  
\roster 
\item
If $\dim D\geq 4$ then:  
$$
\ca(AA) = |\det|^s\otimes \ind_{Q_2}^{H_D}(\ca(A))   
$$
\item
If $\dim D  =2$ then: 
$$
\ca(AA) = |\det|^s\otimes \ind_{Q^+}^{SL_3}(\ca(A^+))+
|\det|^s\otimes\ind_{Q^-}^{SL_3}(\ca(A^-))
$$
\item  
$$
\ca(BB)= |\det|^s\otimes \ind_{QQ\times Q_1}^{GL_2(F)  
\times H_D}(\ca(B))
$$
\item 
$$
V_{N_D}\cong |\det|^t \otimes V_M   
 + |\det|^s \otimes 1 
$$
where $V_M$ is the minimal representation of $M_D$ (center acting 
trivially). 
\endroster 
In the above formulas $\det $ is the usual determinant of 
$2\times 2$ matrices,  
and $s$ and $t$ are given by the following table: 
$$
\matrix
G & s & t \\
E_6 & 2 & 3/2  \\
E_7 & 3 & 2 \\
E_8 & 5 & 3  
\endmatrix
$$
\endproclaim  

{\it Remark:}  In this section we have assumed that $G_D$, and hence
$H_D$ are simply connected. In fact, the isogeny class of $H_D$ 
determines the isogeny class of $G_D$ and vice-versa. 
Note that the action of $H_D$ on 
$\ca(\Om)$ is geometric. It comes from the (algebraic) action of $H_D$ on 
$J_D^0$ which extends uniquely to the group of 
adjoint type. Therefore we have a canonical way of extending the 
representation $V$ to groups of adjoint type. 

\head 8. $\Theta$-lifts from  $G_2(F)$   \endhead

In this section we compute  $\Theta$-lifts of spherical  
tempered representations of  $G_2(F)$ in all three cases. 
In the process we also compute the 
normalizing factors (i.e. coefficients $s$) in Theorems 7.1 
and 7.4. We assume that the residual characteristic of $F$ is
odd.

We study the dual pair $G_2(F)\times F_4(F)$ in a simple group 
$G$ of type $E_8$ first. 
Let $\tilde Y$ be a spherical tempered representation of $G_2(F)$. Then 
$\tilde Y=\Ind_{P_2}^{G_2}(\tilde X)$ where 
$\tilde X$ is a spherical representation of $GL_2(F)$ with
the parameter $(\chi^{-1} |\cdot |^{3/2}, \mu^{-1} |\cdot |^{3/2})$.
As before, $\chi$ and $\mu$ must be unitary characters.
Let $Q_2=L_2V_2$ be the maximal parabolic subgroup of $F_4(F)$ 
stabilizing an amber projective line. 
In Section 6 we have shown that the 
action of $L_2$ on the corresponding amber line gives an 
exact sequence 
$$
1\to SL_3 \to L_2 \to GL_2 \to 1.
$$
One checks that $\rho_{V_2}=|\det|^7$, here $\det$ is the usual 
determinant on $GL_2$. 
Let $X$ be a spherical representation of $GL_2$ with 
the parameter $(\chi |\cdot |^{7/2},
\mu |\cdot |^{7/2}).$ Pull $X$ back to $L_2$. Let 
$$
W=\Ind_{Q_2}^{F_4}X. 
$$ 
Note that $W$ is a unitarizable representation of $F_4(F)$. 
It is quite possible that  $W$ is always irreducible but we 
do not know. 

\proclaim{Theorem 8.1}  
Let $\tilde Y$ be a spherical tempered representation of $G_2(F)$ 
and $W$ the representation of $F_4(F)$ defined above. Then 
$\Theta(\tilde Y)\neq \{\}$ and if $W'\in \Theta(\tilde Y)$ 
then $W'$ is a summand of $W$. 
Let $s\in G_2(\Bbb C)$ be the Satake parameter of $\tilde Y$. Then the 
Satake parameter of $W$ is $\Psi(s\times \rho)$ where 
$\Psi : G_2(\Bbb C)\times SO_3(\Bbb C)\to F_4(\Bbb C)$ and 
$\rho \in SO_3(\Bbb C)$ is the Satake parameter of the 
trivial representation of $SL_2(F)$. 
\endproclaim 

\demo{Proof} 
Let $P=MN$ be 
the maximal parabolic subalgebra of $G$ studied in the last 
two sections. We described an embedding of the dual pair 
$G_2(F)\times F_4(F)$ in $G$ such that $G_2(F)\cap P$ is 
the Heisenberg maximal parabolic subgroup $P_2$. 

Yet another embedding of the dual 
pair $G_2(F)\times F_4(F)$ is given by 
the inclusion of Jordan algebras $J_F \to J$ (use the description of 
the exceptional Lie algebras given in Section 6). In this 
case, $G_2(F)\subset M$ and $F_4(F)\cap P=Q_4=L_4V_4$,  
the Heisenberg maximal parabolic subgroup of $F_4(F)$. 
The Levi component $L_4$ is isomorphic to $GSp_6(F)$. 
Note that the inclusion $GSp_6(F)\to M$ induces an isomorphism 
$GSp_6(F)/Sp_6(F)\cong M/[M,M]$. This is easily seen by 
considering the action of $GSp_6(F)$ and $M$ on $Z$, the 
center of both, $V_4$ and $N$. 

Let $V$ be the minimal representation of $G$. As an 
$G_2(F) \times GSp_6(F)$-module 
$$
V_N \cong V_M \otimes |\det|^3 + 1\otimes |\det|^5 
$$
where $\det$ denotes the usual determinant on $GSp_6(F)$ (Prop. 4.1 
\cite {S}).  
Let $\tilde Y=\Ind_{P_2}^{G_2}(\tilde X)$ 
be a spherical tempered representation of 
$G_2(F)$ and let $Y'$ be the Langlands lift of $\tilde Y$ to 
$PGSp_6(F)$. By Theorem 5.4 $\tilde Y\otimes Y'$ is a quotient  
of $V_M$. It follows from the Frobenius reciprocity 
that $\tilde Y\otimes W'$ is a quotient of $V$ for some $W'$, 
a subquotient of 
$$
\Ind_{Q_4}^{F_4} (Y'\otimes |\det|^3).
$$
On the other hand, $\tilde X\otimes W'$ must be a quotient of 
$V_{U_2}$. By Theorem 7.4 
$$
\ca(AA) = \delta(\det)\otimes \ind_{Q_2}^{H_D}(\ca(A))
$$
for a certain character $\delta$ which we shall now 
determine. 
Although we do not know $\delta$, for a generic choice of 
$\chi$ and $\mu$, $\tilde X\otimes W'$ will be a quotient of
$\ca(AA)$.  
Since $\ca(A)$ is a regular representation of 
$GL_2(F)$ twisted by $\delta$, $\tilde X\otimes W'$ must be a
quotient of 
$$
\tilde X\otimes \Ind_{Q_2}^{F_4}X   
$$
where $X$ is a 
representation of $L_2$ pulled back from a representation 
of $GL_2(F)$ with a parameter $(\chi\delta |\cdot |^{-3/2},
\mu\delta |\cdot |^{-3/2}).$ We get that  
$W'$ is a subquotient of both, $\Ind_{Q_4}^{F_4}(Y'\otimes |\det|^3)$
and $\Ind_{Q_2}^{F_4}X$. This immediately implies that 
$$
\delta(\det)=|\det |^5 .  
$$
The knowledge of $\delta$ implies that
any $\tilde X$ is a quotient of $\ca(AA)$ only and, therefore 
$\Theta$-lifts of $\tilde Y$ are the summands of $W$. 

It remains to check the statement about Satake parameters. 
The $L$-group of $F_4(F)$ is $F_4(\Bbb C)$. Let $Q_3=L_3V_3$ be the 
maximal parabolic subgroup of $F_4$ such that $[L_3,L_3] 
=SL_2 \times SL_3$ where $SL_2$ is spanned by long roots and 
$SL_3$ by short roots. Since $L_3(\Bbb C)$ is the $L$-group of
$L_2$, it fits into exact sequence 
$$
1\to GL_2 \to L_3 \to PGL_3 \to 1.
$$
Let $s\in GL_2(\Bbb C)$ be the parameter $(\chi, \mu)$. Obviously, 
the Satake parameter of $W$ is $s \times \rho \in GL_2 \times SO(3)
\subset L_3$. Let $G_2$ be the centralizer of $SO(3)$ in $F_4$. 
Then $L_3 \cap G_2 = GL_2$ and $s$ is precisely the Satake 
parameter of $\tilde Y$ as in Section 4. 
\enddemo 

Theorem 7.4 can be used, in a similar way, to prove converses of 
Theorems 4.4 and 5.4. We state results without giving details of 
proofs. 

\proclaim{Theorem 8.2} Let $\Phi : SL_3(\Bbb C) \to   G_2(\Bbb C)$ be
the inclusion of the L-groups of $PGL_3(F)$ and $G_2(F)$. 
Let $\tilde Y$ be a tempered spherical representation of
$G_2(F)$.
The Satake parameter of
$\tilde Y$ is $\Phi(s)$ for some $s$, a Satake parameter of a
tempered spherical representation $\tilde Y$ of $G_2(F)$. Note that  
$\Phi(s)=\Phi(s^*)$ where $s^*$ is a Satake parameter of $Y^*$,
the dual of $Y$. Then  
$$
\Theta (\tilde Y) = \{Y, Y^*\}.
$$
\endproclaim 

\proclaim{Theorem 8.3} Let $\Phi : G_2(\Bbb C)\to Spin_7(\Bbb C)$ be
the inclusion of the L-groups of $G_2(F)$ and $PGSp_6(F)$. 
Let $\tilde Y$ be a tempered spherical representation of
$G_2(F)$. Then 
$$
\Theta (\tilde Y) = \{Y' \} 
$$
where $Y'$ is the spherical tempered representation of 
$PGSp_6(F)$ whose Satake parameter is $\Phi(s)$. 
\endproclaim

\head Acknowledgments \endhead

It is a pleasure to thank B. Gross for initiating this project and 
on many discussions during the course of this work. 

A part of this work was written in the Spring of 1995. During that 
time the second author enjoyed hospitality of 
the Mathematical Sciences Research Institute 
and was generously supported by an 
NSF Postdoctoral Fellowship and a Sloan Research Fellowship.

\enddocument 

\bye